\documentclass{gtpart}     
\usepackage{latexsym}
\usepackage{rotating}
\usepackage{amsfonts}
\usepackage{epsfig}
\usepackage{color}
\usepackage{pinlabel}
\input xy
\xyoption{all}

\title{Excision for deformation $K$--theory of free products}

\author{Daniel A. Ramras}
\givenname{Daniel A}
\surname{Ramras}
\address{Department of Mathematics\\\newline
1326 Stevenson Center\\\newline
Vanderbilt University\\\newline
Nashville, TN 37240 U.S.A.}
\email{daniel.a.ramras@vanderbilt.edu}
\urladdr{http://www.math.vanderbilt.edu/~ramrasda}

\arxivreference{math.KT/0703463}
\arxivpassword{q3u4i1l}

\volumenumber{}
\issuenumber{}
\publicationyear{}
\papernumber{}
\startpage{}
\endpage{}
\doi{}
\MR{}
\Zbl{}
\received{}
\revised{}
\accepted{}
\published{}
\publishedonline{}
\proposed{}
\seconded{}
\corresponding{}
\editor{}
\version{}

\newtheorem{theorem}{Theorem}[section]
\newtheorem{lemma}[theorem]{Lemma}
\newtheorem{corollary}{Corollary}[section]

\makeatletter
	\let\c@corollary=\c@theorem
\makeatother

\newtheorem{proposition}{Proposition}[section]

\makeatletter
	\let\c@proposition=\c@theorem
\makeatother

\newtheorem{example}[theorem]{Example}
\newtheorem{definition}[theorem]{Definition}
\newtheorem{remark}[theorem]{Remark}

\DeclareMathOperator*{\colim}{colim}
\DeclareMathOperator*{\tele}{telescope}
\DeclareMathOperator*{\holim}{holim}
\DeclareMathOperator*{\oplusdsp}{\oplus}

\newcommand{\bbZ}{\mathbb{Z}}
\newcommand{\bbN}{\mathbb{N}}
\newcommand{\N}{\mathbb{N}}

\newcommand{\bbC}{\mathbb{C}}
\newcommand{\abs}[1]{\left| #1\right|}

\newcommand{\leqs}{\leqslant}
\newcommand{\geqs}{\geqslant}
\newcommand{\heq}{\simeq}

\newcommand{\maps}{\longrightarrow}

\newcommand{\homeo}{\cong}
\newcommand{\isom}{\cong}
\newcommand{\cross}{\times}
\newcommand{\mcC}{\mathcal{C}}
\newcommand{\mcD}{\mathcal{D}}
\newcommand{\mcR}{\mathcal{R}}

\newcommand{\injects}{\hookrightarrow}

\def\co{\colon\thinspace}
\newcommand{\Rep}{\mathrm{Rep}}

\newcommand{\Hom}{\mathrm{Hom}}

\newcommand{\K}{K_{\mathrm{def}}}

\newcommand{\dom}{\mathrm{dom}}
\newcommand{\codom}{\mathrm{codom}}
\newcommand{\Ob}{\mathrm{Ob}}
\newcommand{\Mor}{\mathrm{Mor}}

\newcommand{\hofib}{\mathrm{hofib}}
\newcommand{\Stab}{\mathrm{Stab}}
\newcommand{\Map}{\mathrm{Map}}
\newcommand{\Vect}{\mathrm{Vect}}

\begin{document}

\begin{abstract}    Associated to a discrete group $G$, one has the topological category of finite dimensional (unitary) $G$--representations and (unitary) isomorphisms.  Block sums provide this category with a permutative structure, and the associated $K$--theory spectrum is Carlsson's deformation $K$--theory $\K(G)$.  The goal of this paper is to examine the behavior of this functor on free products.  Our main theorem shows the square of spectra associated to $G*H$ (considered as an amalgamated product over the trivial group) is homotopy cartesian.  The proof uses a general result regarding group completions of homotopy commutative topological monoids, which may be of some independent interest.

\end{abstract}

\maketitle


\section{Introduction}  Associated to a discrete group $G$ one has the category of finite dimensional representations of $G$ and $G$--equivariant linear isomorphisms.  This category $\mcR(G)$ can be given a natural topology, and block sums of matrices provide it with the structure of a topological permutative category.  The $K$--theory spectrum associated to $\mcR(G)$ is the (general linear) deformation $K$--theory of $G$, first introduced by Carlsson~\cite[Section 4.6]{Carlsson-derived-rep}.  (To be precise, Carlsson's original definition was essentially the singular complex on $\mcR(G)$.  The present definition was introduced by Tyler Lawson~\cite[Chapter 6]{Lawson-thesis}.)  Restricting to unitary representations and unitary isomorphisms, one obtains unitary deformation $K$--theory.  Since the results in this paper are valid in both cases, we will be intentionally ambiguous in our notation, denoting both of these spectra by $\K(G)$.
The zeroth homotopy group $\K^0 (G)$ of $\K(G)$ may be described as follows: the representation spaces themselves form a monoid $\Rep(G)$ under block sum, and $\K^0 (G)$ is the group completion of the monoid $\pi_0 (\Rep(G))$ (Lemma~\ref{reduction2}).  

The excision problem in deformation $K$--theory asks, for an amalgamated product $G*_K H$, whether the natural map of spectra
\begin{equation}\label{ex}
\K (G*_K H) \stackrel{\phi}{\maps} \holim \Big(\K(G) \maps \K(K) \longleftarrow \K(H) \Big)
\end{equation}
is a weak equivalence.  Excision may be thought of as the statement that deformation $K$--theory maps (certain) co-cartesian diagrams of groups to homotopy cartesian diagrams of spectra.
The purpose of this article is to prove that deformation $K$--theory is excisive on free products.  Our main result is:

\medskip
{\bf Theorem \ref{excision-free-prod}}\qua
{\sl 
Let $G$ and $H$ be finitely generated discrete groups.
Then the diagram of spectra 
$$\xymatrix{
       {\K(G * H)} \ar[r] \ar[d]
       & {\K(G)} \ar[d]\\           
       {\K(H)} \ar[r] 
       & {\K(\{1\})}  
                     }
$$ 
is homotopy cartesian.
}
\medskip

This theorem holds for both unitary and general linear deformation $K$--theory.

Information about the connectivity of the map $\phi$ provides a relationship between the deformation $K$--groups of $G*_K H$ and the deformation $K$--groups of the factors, via the long exact Mayer--Vietoris sequence in homotopy associated to the homotopy pullback square.  
To illustrate the usefulness of this sequence, we apply our excision result to compute the deformation $K$--theory of $PSL_2(\bbZ) \isom \bbZ/2 * \bbZ/3$ (Proposition~\ref{PSL}).

Another motivation for studying excision in deformation $K$--theory is that excision results can improve our understanding of Atiyah--Segal phenomena.  Classically, the Atiyah--Segal theorem~\cite{Atiyah-char-coh, Atiyah-Segal} states that the representation ring of a compact Lie group, after completion at the augmentation ideal, is isomorphic to the complex $K$--theory ring $K^*(BG)$.  When $G$ is an infinite discrete group, one hopes to prove analogous theorems relating $\K^*(G)$ and $K^*(BG)$ (here we are thinking of unitary deformation $K$--theory).  This has been done in certain cases.  For free groups, T. Lawson~\cite[p. 11]{Lawson-simul} has shown that $\K(G)$ is weakly equivalent, as a $\mathbf{ku}$--module, to $\mathbf{ku}\vee \Sigma \mathbf{ku}$.  For fundamental groups of compact, aspherical surfaces $M$, the author has established an isomorphism in homotopy $\K^*(\pi_1 M) \isom K^*(M)$, for $*>0$~\cite{Ramras-surface}.  (In fact, this isomorphism holds in dimension zero when $M$ is non-orientable.)  The failure in low degrees is an important feature of deformation $K$--theory: while topological $K$--theory $\K^*(BG)$ depends only on the stable homotopy type of $BG$, deformation $K$--theory is a subtler invariant of $G$, closely tied to the topology of the representation spaces $\Hom(G, U(n))$. 

We note that excision is not satisfied for all amalgamated products.  The fundamental group of a Riemann surface, described as an amalgamated product via a connected sum decomposition of the surface, fails to satisfy excision on $\pi_0$ (we expect that the natural map $\phi$ induces an isomorphism in positive degrees, though; at least the domain and range of $\phi$ have abstractly isomorphic homotopy groups, by the main result of~\cite{Ramras-surface}).  For further discussion, we refer the reader to~\cite[Section 6]{Ramras-surface}.  

We briefly describe our approach to the excision problem.  Since deformation $K$--theory arises from a permutative category $\mcR(G)$, its zeroth space is (weakly equivalent to) the group completion
$\Omega B\abs{\mcR(G)}$ of the monoid $\abs{\mcR(G)}$.
McDuff and Segal, in their well-known paper~\cite{McDuff-Segal}, introduced a simple model for the homology type of the group completion.  This model arises as an infinite mapping telescope $M_{\infty} (m)$ of the monoid $M$, in which the maps come from multiplication by a fixed element $m\in M$.  Ordinarily, starting from this homological model one tries to apply Quillen's plus-construction to abelianize the fundamental group and produce a space weakly equivalent to $\Omega BM$.  We show  that under certain simple conditions, the fundamental group of this mapping telescope is actually abelian, and as a consequence this space is weakly equivalent to $\Omega BM$.  Precisely, our main technical result is as follows.

\medskip
{\bf Theorem \ref{gp-completion-thm}}\qua
{\sl 
Let $M$ be a homotopy commutative monoid which is stably group-like with respect to an anchored element $m\in M$.  Then there is a natural isomorphism
$$\eta\co \pi_* M_{\infty} (m)\stackrel{\isom}{\maps} \pi_* \Omega BM,$$
and the induced map on $\pi_0$ is an isomorphism of groups.  The map $\eta$ is induced by a zig-zag of natural weak equivalences.
}
\medskip

The term \emph{anchored} is introduced in Definition~\ref{anchored}.
Theorem~\ref{excision-free-prod} is an application of this general result, and this result also forms the basis for the computations in~\cite{Ramras-surface} of $\K^*(\pi_1 (\Sigma))$ and $\pi_1 \Hom(\pi_1 \Sigma, U)/U$ for $\Sigma$ a compact, aspherical surface.

The paper is organized as follows.  In Section~\ref{K-def}, we introduce deformation $K$--theory, describe its zeroth space, and prove a result of T. Lawson~\cite[Chapter 6.2]{Lawson-thesis} relating its underlying monoid to homotopy orbits of representation spaces.  In Section~\ref{gp-completion}, we study group completions of topological monoids, and prove our results regarding the McDuff--Segal model for the homology of $\Omega BM$.  These results are applied to deformation $K$--theory in Section~\ref{gp-completion-2}.  The main result regarding excision for free products appears in Section~\ref{free-excision}, as does the computation of $\K^*(PSL_2 (\bbZ))$.

\vspace{.15in}

\noindent {\bf Acknowledgements:}  

This work formed part of the author's Stanford University Ph.D. thesis, and most of this paper was written while the author was a graduate student at Stanford.  I would like to thank my advisor, Gunnar Carlsson, for suggesting this topic and for all his help and encouragement.  I would also like to thank Tyler Lawson and Soren Galatius for helpful conversations, and the referee, whose comments greatly improved the exposition.  Much of this paper was written during visits to Columbia University, and I would like to thank the Columbia Mathematics Department for its hospitality.

This research was supported by an NSF Graduate Research Fellowship.


\section{Deformation $K$--theory}$\label{K-def}$

In this section, we introduce Carlsson's deformation $K$--theory and discuss its basic properties.  Deformation $K$--theory is a contravariant functor from discrete groups to spectra, and is meant to capture homotopy theoretical information about the representation varieties of the group in question.  We will construct a connective $\Omega$--spectrum $\K(G)$ by considering the $K$--theory of an appropriate permutative topological category of representations (this category was first introduced by Lawson~\cite[Chapter 6]{Lawson-thesis}).  Although we phrase everything in terms of the unitary groups $U(n)$, all of the constructions, definitions and results in this section are valid for the general linear groups $GL_n(\bbC)$, and only notational changes are needed in the proofs.

For the rest of this section, we fix a discrete group $G$.

\begin{definition}  Associated to $G$ we have a topological category (i.e. a category object in the category of compactly generated spaces) $\mcR(G)$ with object space
$$\Ob (\mcR(G)) = \coprod_{n=0}^{\infty} \Hom(G, U(n))$$
and morphism space
$$\Mor (\mcR(G)) = \coprod_{n=0}^{\infty} U(n) \cross \Hom(G,U(n)).$$
\end{definition}

The domain and codomain maps are 
$\dom(A, \rho) = \rho$ and 
$\codom(A, \rho) = A \rho A^{-1}$, 
and composition is given by $(B, A\rho A^{-1})\circ (A, \rho) = (BA, \rho)$. 
The identity map $\Ob (\mcR(G)) \to \Mor (\mcR(G))$ sends a representation $\rho\in \Hom(G, U(n))$
to the pair $(\rho, I_n)$, where $I_n\in U(n)$ denotes the identity matrix.
The representation spaces are topologized using the compact--open topology, or equivalently as subspaces of $\prod_{g\in G} U(n)$.  
We define $U(0)$ to be the trivial group, and the single point $* \in \Hom(G, U(0))$ will serve as the basepoint.

The functor
$\oplus\co \mcR(G)\cross \mcR(G)\to \mcR(G)$
defined via block sums of unitary matrices is continuous and strictly associative, with the trivial representation $*\in \Hom(G, U(0))$ as unit, and this functor makes $\mcR(G)$ into a permutative category in the sense of~\cite{May-perm}.  Explicitly, given representations $\rho\in \Hom(G, U(n))$ and $\psi\in \Hom(G, U(m))$, the representation $\oplus(\rho, \psi) = \rho\oplus \psi \in \Hom(G, U(n+m))$ is defined, for each $g\in G$, by
$$\rho\oplus \psi (g) = \left[ 
\begin{array}{cc}
       \rho (g) & 0 \\           
       0 & \psi (g) 
\end{array}
\right].
$$
The behavior of $\oplus$ on the morphisms of $\mcR(G)$ is defined similarly.  
We note that the single point in $\Hom(G, U(0))$ acts as the identity under the operation $\oplus$.
The natural commutativity isomorphism
$$c\co \rho \oplus \psi  \stackrel{\isom}{\maps} \psi \oplus \rho$$
is defined via the (unique) permutation matrices $\tau_{n,m}$ satisfying 
$$\tau_{n,m} (A\oplus B) \tau_{n,m}^{-1} = B\oplus A$$
for all $A\in U(n)$ and $B\in U(m)$.  (A general discussion of the functor associated to a collection of matrices like this one can be found in the proof of Corollary~\ref{gp-completion-cor}.)
Each homomorphism $f\co G\to H$ induces a continuous functor $f^*\co \mcR(H)\to \mcR(G)$, 
defined on objects by $\rho \mapsto \rho \circ f$ and on morphisms by $(\rho, A) \mapsto (\rho\circ f, A)$.
It follows from the definitions that this functor is permutative, and that $\K(-)$ defines a contravariant functor from the
category of discrete groups and group homomorphisms to the category of topological permutative categories and continuous, permutative functors.

May's machine~\cite{May-perm} constructs  a (special) $\Gamma$--category (in the sense of~\cite{Segal-cat-coh}) associated to any permutative (topological) category $\mcC$.  (We note that May's construction requires the identity object in $\mcC$ to be a non-degenerate basepoint; in our case the identity object is the unique zero-dimensional representation and is \emph{disjoint} from the rest of $\Ob (\mcR(G))$.)
Taking geometric realizations yields a special $\Gamma$--space, and Segal's machine then produces a connective $\Omega$--spectrum $K(\mcC)$, the $K$--theory of the permutative category $\mcC$.  This entire process is functorial in the permutative category $\mcC$.

\begin{definition}$\label{K-def-spectrum}$
Given a discrete group $G$, the deformation $K$--theory spectrum of $G$ is defined to be the 
$K$--theory spectrum of the permutative category $\mcR(G)$, i.e.
$$\K(G) = K(\mcR(G)).$$
This spectrum is contravariantly functorial in $G$.
\end{definition}

We now describe the zeroth space of the spectrum $\K(G)$.  

\begin{lemma}$\label{zeroth-space}$  For each discrete group $G$, the zeroth space of $\K (G)$ is naturally weakly equivalent to $\Omega B (\abs{\mcR(G)})$, where $B$ denotes the bar construction on the topological monoid $\abs{\mcR(G)}$.
\end{lemma}

The proof of this result is just an elaboration of the May's proof~\cite[Construction 10, Step 2]{May-perm}  that the $\Gamma$--category associated to a permutative category $\mcC$ is special.  

\vspace{.15in} Next, we discuss an observation due to Lawson~\cite[Chapter 6.2]{Lawson-thesis} regarding the classifying space of the category $\mcR(G)$.  For convenience of the reader, and to set notation, we include a discussion of the simplicial constructions of the classifying space $BU(n)$ and the universal bundle $EU(n)$.  This discussion is based on Segal~\cite[Section 3]{Segal-class-ss}.

Associated to $G$ we have the homotopy orbit spaces
$$\Hom(G, U(n))_{hU(n)} = EU(n)\cross_{U(n)} \Hom(G, U(n)),$$
where $EU(n)$ denotes the total space of a universal (right) principal $U(n)$--bundle,
and $U(n)$ acts on $EU(n) \cross \Hom(G, U(n))$ via $(e,\rho)\cdot g = (eg, g^{-1}\rho g)$.  
We want to define a monoid structure on the disjoint union (over $n\in \bbN$) of these spaces.

We take $EU(n)$ to be the classifying space of the translation category $\overline{U(n)}$ of $U(n)$, that is, the topological category whose object space is $U(n)$ and whose morphism space is $U(n)\cross U(n)$.  The morphism $(A, B)$ is the unique morphism from $B$ to $A$ in $\overline{U(n)}$, and hence
composition is given by $(A,B)(B,C) = (A,C)$.  We will write $(A,B) = B\to A$.

Since morphisms in $\overline{U(n)}$ are determined by their domain and codomain, functors into $\overline{U(n)}$ are determined by their behavior on objects.  For each $g\in U(n)$, we have a functor
$$r_g \co \overline{U(n)} \maps \overline{U(n)}$$
defined on objects by $r_g (A) = Ag$.  Since $r_{gh} = r_h \circ r_g$, these functors define a \emph{right}
action of $U(n)$ on $\overline{U(n)}$.  After geometric realization, this defines a continuous right action of $U(n)$ on $EU(n) = \abs{\overline{U(n)}}$.  

The quotient map $EU(n) \to EU(n)/U(n)$ associated to this action is a universal, principal $U(n)$--bundle (Segal~\cite[Section 3]{Segal-class-ss}).  Since the action of $U(n)$ on $EU(n)$ arises from a simplicial action of $U(n)$ on the simplicial space $N_\cdot \overline{U(n)}$, the quotient space $EU(n)/U(n)$ may be described as the realization of the simplicial space
$$k\mapsto \left(N_k \overline{U(n)}\right)/U(n).$$
This simplicial space may
also be described categorically, as follows.  Let $\mcC_{U(n)}$ be the topological category with a unique object $*$ and with $U(n)$ as its morphism space, where composition is given by 
$A\circ B = AB$.  We define $BU(n):= \abs{\mcC_{U(n)}}$.
Consider the functor 
$$\overline{U(n)} \stackrel{\pi}{\maps} \mcC_{U(n)}$$
defined on objects by $A\mapsto *$ and on morphisms by $(A, B) \mapsto AB^{-1}$.  The induced
map $N_\cdot \pi$ on nerves factors through the $U(n)$--action on $N_k \overline{U(n)}$, inducing
a map of simplicial spaces
$$\left(N_k \overline{U(n)} \right) /U(n) \maps N_k \mcC_{U(n)}$$
which is a homeomorphism for each $k\in \bbN$.  
Hence the induced map on realizations
$$EU(n)/U(n) \maps \abs{\mcC_{U(n)}} = BU(n)$$
is a homeomorphism.  Thus $EU(n) \stackrel{\abs{\pi}}{\maps} BU(n)$
is a universal principal $U(n)$--bundle.

We now turn to the monoid structure on the space
$$\coprod_{n=0}^{\infty} \Hom(G, U(n))_{hU(n)}
:= \coprod_{n=0}^{\infty} EU(n)\cross_{U(n)} \Hom(G, U(n)),$$
which we (abusively) denote by $\Rep(G)_{hU}$.
The continuous block sum maps
$\oplus\co U(n)\cross U(m) \to U(n+m)$
yield continuous functors 
$$\oplus\co \overline{U(n)}\cross \overline{U(m)}\maps \overline{U(n+m)}$$
defined on objects by
$\oplus (A, B) = A\oplus B$.  By abuse of notation, we will denote the realizations of these functors
by $\oplus\co EU(n)\cross EU(m)\to EU(n+m)$.
The multiplication on $\Rep(G)_{hU}$ will also be denoted by $\oplus$, and is defined as follows.
Each point in this monoid is represented some pair $(e, \rho)$, where $e\in EU(n)$ and 
$\rho\in \Hom(G, U(n))$ (for some $n\in \bbN$).  We set 
$$[e_1, \rho_1] \oplus [e_2, \rho_2] = [e_1\oplus e_2, \rho_1 \oplus \rho_2].$$
It follows from the definitions that this map is well-defined and continuous, and yields an associative multiplication on $\Rep(G)_{hU}$ with the single point in dimension zero as the identity.

We may now define a homomorphism of topological monoids
$$\alpha \co \Rep(G)_{hU} \stackrel{\alpha}{\maps} \abs{\mcR(G)}$$
as follows.  Each point in $EU(n)$ is represented by a pair $(x, t)\in N_k \overline{U(n)} \cross \Delta^k$
(for some $k\in \bbN$), where $\Delta^k$ denotes the topological simplex of dimension $k$.
Since $N_k \overline{U(n)}\isom U(n)^{k+1}$, we may write
$x = A_0 \to A_1\to \cdots \to A_k$ for some $A_i\in U(n)$.  Hence each point in $\Rep(G)_{hU}$ is represented
by a triple $(A_0 \to A_1\to \cdots \to A_k, t, \rho)$ with $A_i \in U(n)$, $t\in \Delta^k$, and $\rho\in \Hom(G, U(n))$
(for some $n, k \in \bbN$).  Recall that morphisms in $\mcR(G)$ are determined by their domain (a representation $\rho$) and a unitary matrix $A$.  We write such a morphism as
$$\rho \stackrel{A}{\maps} A\rho A^{-1}.$$
We now define
\begin{equation*}\begin{split}
\alpha & \left(  [A_0 \to A_1\to \cdots \to A_k, t, \rho] \right) \\
&= [A_0 \rho A_0^{-1} \xrightarrow{A_1 A_0^{-1}} A_1 \rho A_1^{-1} 
\xrightarrow{A_2 A_1^{-1}} \cdots \xrightarrow{A_k A_{k-1}^{-1}} A_k \rho A_k^{-1}, t].
\end{split}
\end{equation*}								
By tracing the definitions, one may check that this map is a well-defined, continuous homomorphism of monoids.

Lawson's observation, then, is:

\begin{proposition}[Lawson]$\label{classifying-space}$
The map $\alpha\co Rep(G)_{hU} \maps  \abs{\mcR(G)}$ is an isomorphism of topological monoids.
\end{proposition}
\begin{proof}  We will factor $\alpha$ as a composition
of homeomorphisms. 

We begin by considering $\Hom(G, U(n))$ as a constant simplicial space, so that
$$EU(n)\cross \Hom(G, U(n)) \homeo \left| k\mapsto N_k \overline{U(n)} \cross \Hom(G, U(n)) \right|.$$
Now, combining the level-wise action of $U(n)$ on $EU(n)$ with the conjugation action of $U(n)$ on $\Hom(G, U(n))$ gives the simplicial space on the right a simplicial right action of $U(n)$, and we have a homeomorphism
$$EU(n)\cross_{U(n)} \Hom(G, U(n)) 
	\homeo \left| k\mapsto \left( N_k \overline{U(n)} \cross \Hom(G, U(n)) \right)/U(n) \right|.
$$
The continuous functions
$$f_k \co \left( N_k \overline{U(n)} \cross \Hom(G, U(n)) \right)/U(n) \maps
N_k \mcR(G)$$
defined by
\begin{equation*}\begin{split}
f_k & \left(  [A_0 \to A_1\to \cdots \to A_k, \rho] \right) \\
&= A_0 \rho A_0^{-1} \xrightarrow{A_1 A_0^{-1}} A_1 \rho A_1^{-1} 
\xrightarrow{A_2 A_1^{-1}} \cdots \xrightarrow{A_k A_{k-1}^{-1}} A_k \rho A_k^{-1}
\end{split}
\end{equation*}combine to form a simplicial map, which we claim is a homeomorphism on realizations.
In fact,  each $f_k$ is bijective, with inverse given by
\begin{equation*}\begin{split}
f_k^{-1} \left( \rho \stackrel{A_1}{\maps} A_1 \rho A_1^{-1} \right. & 
\left. \stackrel{A_2}{\maps} A_2 A_1 \rho A_1^{-1} A_2^{-1} 
\stackrel{A_3}{\maps} \cdots \stackrel{A_k}{\maps} A_k \cdots A_1 \rho A_1^{-1} \cdots A_k ^{-1}\right)\\
&  = [I \to A_1 \to A_2 A_1 \to A_k \cdots A_1, \rho].
\end{split}
\end{equation*}
Since the domain of $f_k$ is compact and the range is Hausdorff, we conclude that each $f_k$ is a homeomorphism, and consequently the induced map on realizations is a homeomorphism as well.
\end{proof}

We end this section with a simple observation regarding the zeroth homotopy group of deformation $K$--theory.  The topological monoid $\Rep(G)$ is defined by
$$\Rep(G) = \coprod_{n=0}^{\infty} \Hom(G, U(n)).$$ 
The monoid structure on $\Rep(G)$ is given by block sum of representations, just as in the definition of the permutative structure on the category $\mcR(G)$.  (In fact, $\Rep(G)$ is precisely the submonoid of 
$\abs{\mcR(G)}$ corresponding to the subcategory of identity morphisms.)
Again, $U(0)$ is the trivial group and the single element in $\Hom(G, U(0))$ will act as the identity.  The same construction may be applied with the general linear groups in place of the unitary groups, and we keep the notation intentionally vague.  

\begin{lemma}$\label{reduction2}$
Let $G$ be a discrete group.  Then $\K^0 (G) \homeo Gr(\pi_0 (\Rep(G)))$, where $Gr$ denotes the group-completion of a monoid, i.e. its Grothendieck group.
\end{lemma}
\begin{proof}  By Lemmas~\ref{zeroth-space} and~\ref{classifying-space}, we know that $\K^0(G)$ is the group completion of the monoid
$\pi_0 (\Rep(G))_{hU}$, so we just need to show that there is an isomorphism of monoids
$\pi_0(\Rep(G))_{hU}\isom \pi_0(\Rep(G))$.
This follows from the fact that the monoid maps
$$\Rep(G)_{hU} \longleftarrow \coprod_{n=0}^{\infty} EU(n)\cross \Hom(G, U(n)) \maps \Rep(G)$$
are both fibrations, with connected fibers.
\end{proof}


\section{Group completion}$\label{gp-completion}$

The goal of this section is to provide a convenient homotopy theoretical model for the group completion of a topological monoid satisfying certain simple properties.  The results of this section will be applicable to deformation $K$--theory, and form the basis of our excision results as well as the computations in Ramras~\cite{Ramras-surface}.

The models for group completion that we will study arise as mapping telescopes, as in McDuff--Segal~\cite{McDuff-Segal}.  Throughout this section $M$ will denote a homotopy commutative topological monoid and $e\in M$ will denote the identity element.  We write the multiplication in $M$ as $*$, and for each $m\in M$ we denote the $n$--fold product of $m$ with itself by $m^n$.

\begin{definition}  For each $m\in M$, we denote the mapping telescope 
$$\tele (\underbrace{M\stackrel{* m}{\maps} M\stackrel{* m}{\maps} 
                         \cdots \stackrel{* m}{\maps} M}_N)$$
by $M_N (m)$, and we denote the infinite mapping telescope 
$$\colim_{N\to \infty} M_N (m) = \tele (M\stackrel{* m}{\maps} M\stackrel{* m}{\maps} \cdots)$$
by $M_{\infty} (m)$.
\end{definition}

We denote points in these telescopes by triples $(x, n, t)$, where $x\in M$, $n\in \bbN$ and $t\in [0,1)$.  Note that each of these spaces is functorial in the pair $(M, m)$, in the sense that if
$M\to M'$ is a map of monoids, then there are induced maps of telescopes
$$M_N (m) \maps M'_N(f(m))$$
(for each $N\in \bbN\cup \{\infty\}$).
To be precise, there is a category $\mathcal{M_*}$ whose objects are
pairs $(M, m)$, where $M$ is a monoid and $m\in M$ is any element, and whose morphisms 
$(M, m) \to (M', m')$ are homomorphisms $f\co M\to M'$ with $f(m) = m'$; for every $N\in \bbN\cup \{\infty\}$, the assignment 
$$(M, m) \mapsto \left( M_N(m), (e, 0, 0) \right)$$ 
extends to a functor from $\mathcal{M}_*$ to the category of \emph{based} spaces. 
We will denote the basepoint $(e,0,0)\in M_N (m)$ simply by $e$.

\begin{definition} We say that $M$ is \emph{stably group-like} with respect to an element $m\in M$ if the cyclic submonoid of $\pi_0(M)$ generated by $m$ is cofinal.  In other words, $M$ is stably group-like with respect to $m$ if for every $x\in M$ there exists $y\in M$ and $n\in \bbN$ such that $x* y$ and $m^n$ lie in the same path component of $M$.  We refer to such $y$ as \emph{stable homotopy inverses} for $x$ (with respect to $m$).
\end{definition}

The reason for our terminology is the following result.

\begin{proposition}$\label{stbly-gp-like}$
Let $(M, *)$ be a homotopy commutative monoid and let $m\in M$ be an element.  Then there is a natural abelian monoid structure on $\pi_0 (M_{\infty}(m))$, and $M$ is stably group-like with respect to $m$ if and only if $\pi_0 (M_{\infty} (m))$ is a group under this multiplication.  If $M$ is stably group-like with respect to $m$, then $\pi_0 (M_{\infty})$ is the group completion of $\pi_0 (M)$.
\end{proposition}
\begin{proof} Given components $C_1$ and $C_2$ in $\pi_0 (M_{\infty} (m))$, choose representatives $(x_1, n_1, 0)$ and $(x_2, n_2, 0)$ for $C_1$ and $C_2$ respectively.  We define
$C_1* C_2$ to be the component containing $(x_1* x_2, n_1+n_2, 0)$.  To see that this operation is well-defined, one uses the fact that if $(x,n,0)$ and $(x',n',0)$ are connected by a path, then this path lies in some finite telescope $M_N (m)$ (with $N>n, n'$) and one can collapse the first $N$ stages of the telescope to obtain a path in $M$ from $x* m^{N-n}$ to $x'* m^{N-n'}$.  
Since this proposition will not be needed in the sequel, we leave the rest of the proof as an exercise for the reader.
\end{proof}

\begin{example}$\label{fin-gen}$
Let $(M, *)$ be a homotopy commutative monoid and assume that $\pi_0 (M)$ is finitely generated, with generators $m_1, \ldots, m_k\in M$.  Then $M$ is stably group-like with respect to the element $m = m_1* \ldots * m_k$: each component is represented by a word in the $m_i$, and we may add another word to even out the powers.  
More precisely, if $c = m_{i_1}^{a_1} * \cdots * m_{i_k}^{a_k}$, then let
$A = \max_j a_j$.  The element $c^{-1} = m_{i_1}^{A-a_1} * \cdots * m_{i_k}^{A - a_k}$ is a stable homotopy inverse for $c$ with respect to $m$. 
(This example appears, in spirit at least, in McDuff and Segal~\cite[p. 281]{McDuff-Segal}, and will be central to our results on excision.)
\end{example}

Before stating the main result of this section, we need the following definition.

\begin{definition}$\label{anchored}$ Let $(M, *)$ be a homotopy commutative monoid.  We call an element $m\in M$ 
\emph{anchored} if there exists a homotopy $H\co M\cross M\cross I\to M$ such that for every $m_1, m_2\in M$, 
$H_0(m_1, m_2) = m_1* m_2$, $H_1(m_1, m_2) = m_2* m_1$, and $H_t(m^n, m^n) = m^{2n}$ for all $t\in I$ and all $n\in \N$. 

We define the subcategory $\mathcal{M}_*^{\mathrm{asg}} \subset \mathcal{M}_*$ to be the full subcategory on the objects $(M, m)$ such that $m$ is anchored in $M$ and $M$ is stably group-like with respect to $m$.
\end{definition}

\begin{theorem}$\label{gp-completion-thm}$ Let $M$ be a homotopy commutative monoid which is stably group-like with respect to an anchored element $m\in M$.  Then there is a natural isomorphism
$$\eta\co \pi_* M_{\infty} (m)\stackrel{\isom}{\maps} \pi_* \Omega BM,$$
and the induced map on $\pi_0$ is an isomorphism of groups.  The map $\eta$ is induced by a zig-zag of natural weak equivalences between functors from $\mathcal{M}_*^{\mathrm{asg}}$ to based spaces.
\end{theorem}

We note that the final statement of the theorem is slightly stronger than saying that $M_{\infty} (m)$
and $\Omega BM$ are naturally isomorphic in the homotopy category: the naturality diagrams for this zig-zag are strictly commutative, not just homotopy commutative.  This strict commutativity will be used in the proof of our main result, Theorem~\ref{excision-free-prod}.

A number of comments are in order regarding zig-zags, basepoints, and the precise meaning of naturality for the map $\eta$.  By a zig-zag we simply mean a sequence of spaces $X_1, X_2, \ldots, X_k$, together with maps $f_i$ between $X_i$ and $X_{i+1}$ (in either direction).  The sequence of spaces in our natural zig-zag is described in (\ref{gc-zig-zag}) below, and involves only one auxiliary functor, denoted by $\hofib (q_m)$.  A map $f\co X\to Y$ between possibly disconnected spaces will be called a weak equivalence if and only if it induces isomorphisms $f_*\co \pi_* (X, x) \to \pi_* (Y, f(x))$ for all $x\in X$.  

The isomorphism $\eta$ on homotopy groups will be valid for all compatible choices of basepoint, in the following sense.  The zig-zag of isomorphisms on $\pi_0$ gives an isomorphism
$\eta_0\co \pi_0 M_{\infty} (m) \to \pi_0 \Omega BM$, and we call basepoints $x\in M_{\infty}(m)$ and $y\in \Omega BM$ compatible if
$\eta_0 ([x]) = [y]$.  Now, for each pair $x$ and $y$ of compatible basepoints there is in fact a canonical isomorphism
$$\eta_{x,y}\co  \pi_* \left (M_{\infty} (m), x\right) \stackrel{\isom}{\maps} \pi_* \left( \Omega BM, y \right).$$
This isomorphism is constructed using the fact that if $X$ is a simple space, meaning that the action of $\pi_1(X,x)$ on $\pi_n (X, x)$ is trivial for every $n\geqs 1$, then any two paths between points $x_1, x_2\in X$ induce the same isomorphism $\pi_* (X, x_1)\to \pi_* (X, x_2)$.   Since we are dealing with a zig-zag of weak equivalences ending with a simple space, all spaces involved are simple, and hence $\eta_{x,y}$ is well-defined.

Naturality for the map $\eta$ means that for each morphism $f\co(M, m) \to (N, f(m))$ in $\mathcal{M}^{\mathrm{asg}}_*$, and for each pair of compatible basepoints $x\in M_{\infty} (m)$ and
$y\in \Omega BM$, we have 
$$\Omega Bf \circ \eta_{x,y} = f_{\infty} \circ \eta_{f_{\infty} (x), (\Omega Bf) (y)}\co 
	\pi_* \left( M_{\infty} (m), x\right) \maps \pi_*\left( \Omega BN, (\Omega Bf)(y)\right),$$
where $f_{\infty} = f_{\infty} (m)$ denotes the map on telescopes induced by $f$.
This equation follows from the naturality of the weak equivalences involved in the zig-zag (\ref{gc-zig-zag}).

\begin{remark}
It is possible to relax the definition of ``anchored" without affecting Theorem~\ref{gp-completion-thm} (and only minor changes are needed in the proof).  For example, the homotopies anchoring $m^n$ need not be the same for all $n$, and in fact we only need to assume their existence for ``enough'' $n$.  (In particular, it is not necessary to assume that there is a homotopy anchoring $m^0 = e$.)  For all our applications, though, the current definition suffices.
\end{remark}

We now turn to the proof of Theorem~\ref{gp-completion-thm}.  To fix notation, we begin by describing the McDuff--Segal approach to group completion~\cite{McDuff-Segal}.
Given a space $X$ together with an (left) action of a monoid $M$ on $X$, one may form the topological category $X_M$ whose object space is $X$ and whose morphism space is $M\cross X$.  Here $(m,x)$ is a morphism from $x$ to $m\cdot x$, and composition is given by $(n,mx)\circ (m,x) = (n*  m, x)$.
There is a natural, continuous functor $Q\co X_M\to BM$ where $BM$ denotes the topological category with one object and with morphism space $M$ (the geometric realization of $BM$ is the classifying space of $M$, which we also denote by $BM$).  On morphisms, this functor sends $(m, x)$ to $m$.
When $X = M$ (acted on via left multiplication) the category $M_M$ has an initial object (the identity $e\in M$) and hence $EM = \abs{M_M}$ is canonically contractible.  (Note here that $M_M$ is \emph{not} the category with a unique morphism between each pair of objects).  Now, $M$ acts on $M_{\infty} (m)$ via $x\cdot (y,n,t) = (x* y, n, t)$, and we 
define $(M_{\infty})_M = (M_{\infty} (m))_M$.  This space has a natural basepoint, coming from the basepoint $e\in M_{\infty} (m)$.  Observe that $(M_{\infty})_M$ is the infinite mapping telescope of 
the sequence
$$EM\stackrel{F_m}{\maps} EM \stackrel{F_m}{\maps}  \cdots ,$$
where $F_m$ is the functor  defined by $F_m (x) = x* m$ and  $F_m (n,x) = (n, x* m)$.  Since $EM$ is contractible, it follows that $(M_{\infty})_M$ is (weakly) contractible as well.  Now, as noted above we have a functor $Q_m\co (M_{\infty})_M\to BM$, and we denote its realization by $q_m$.  The fiber of the map $q_m$ (over the vertex of $BM$) is precisely $M_{\infty} (m)$, and so we have a natural map
$$i_m\co M_{\infty} (m) \maps \hofib(q_M).$$

The theorem of McDuff and Segal~\cite[Proposition 2]{McDuff-Segal} states that 
this map induces an isomorphism 
$$(i_m)_*\co H_* (M_{\infty} (m), (i_m)^*(A)) \stackrel{\isom}{\maps} H_* (\hofib (q_m), A)$$
in homology with local (abelian) coefficients, so long as the action of $M$ on $M_{\infty} (m)$ is by homology equivalences (again with abelian local coefficients).  
This hypothesis is satisfied when $M$ is stably group-like with respect to $m$;
this is essentially an exercise in the definitions, using the fact that homology of a mapping telescope may be computed as a colimit (for a complete proof, see Ramras~\cite[Lemma 3.0.25]{Ramras-thesis}).  

Next, since $(M_{\infty})_M$ is (weakly) contractible, we have a weak equivalence from $\Omega BM$ to $\hofib(q_M)$, induced by the diagram
$$\xymatrix{
       {*} \ar[r]^(.35){\heq} \ar[d] &  (M_{\infty})_M \ar[d]^{q_m}\\
       BM \ar@{=}[r]   & {BM.}  
                     }
$$ 
Here the maps from $*$ are the inclusions of the natural basepoints; note that
$\Omega BM \homeo \hofib(*\to BM)$.  Hence we have a zig-zag of natural transformations
\begin{equation}\label{gc-zig-zag}
M_{\infty} (m) \stackrel{i_m}{\maps} \hofib (q_m) \stackrel{\heq}{\longleftarrow} \Omega BM,
\end{equation}
and the first map induces an isomorphism in homology with local (abelian) coefficients.  (This is the full conclusion of the McDuff--Segal Theorem.)

Ordinarily, one would now attempt to show that after applying a plus-construction to the space $M_{\infty} (m)$, this space becomes weakly equivalent to the other two.  We will show, though, that when $m$ is anchored in $M$, the fundamental group $\pi_1 (M_{\infty} (m), x)$ is already abelian, and hence no plus-construction is required (this holds for every basepoint $x$).  This will allow us to deduce Theorem~\ref{gp-completion-thm}. 

\begin{remark}
We note that McDuff and Segal actually work with the thick realization $||\cdot ||$ of simplicial spaces, meaning that the real conclusion of their theorem is that the map $M_{\infty} (m) \to \hofib ||Q_m||$ is a homology equivalence with local coefficients.  (Note that $M_{\infty} (m)$ is the fiber of both $||Q_m||$ and $q_m = \abs{Q_m}$.)
In the present application, the simplicial spaces involved are good, so the thick realization is homotopy equivalent to the ordinary realization by Segal~\cite[Proposition A.1]{Segal-cat-coh}.  Hence one finds that there is a weak equivalence $\hofib ||q|| \to \hofib (q)$,
and since weak equivalences induce isomorphisms in homology with local coefficients, we conclude that the map $i_M\co M_{\infty} (m)\to \hofib (q)$ induces isomorphisms in homology with local coefficients as well.  
\end{remark}

We will now show that all components of 
$M_{\infty} (m)$ have abelian fundamental group, and we begin with the component containing $e = (e,0,0)$.

\begin{proposition}$\label{abelian}$
Let $M$ be a homotopy commutative monoid.  If $m\in M$ is anchored, then 
$\pi_1 (M_{\infty}(m), e)$ is abelian.
\end{proposition}

The idea of the proof is show that the ordinary multiplication in $\pi_1 (M_{\infty} (m), e)$ agrees with an operation defined in terms of the multiplication in $M$.  This latter operation will immediately be commutative, by our assumptions on $M$.

The proof will require some simple lemmas regarding loops in mapping telescopes.  We write $p_1 \centerdot p_2$ for composition of paths (tracing out $p_1$ first).
We begin by describing the type of loops that we will need to use.  (Our notation for mapping telescopes was described at the start of this section.)

\begin{figure}[ht!]
\labellist
\small\hair 2pt
\pinlabel $\gamma$ at 65 131
\pinlabel $N(\gamma)$ at 210 59
\endlabellist
\centering
\includegraphics[scale=0.80]{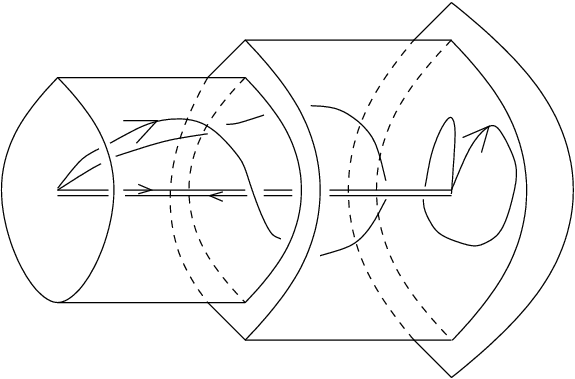}
\caption{\small{A loop $\gamma$ and its normalization $N(\gamma)$.}}
\end{figure}

\begin{definition}$\label{loops}$
For each $n$, there is a canonical path $\gamma_n\co I \to M_{\infty}(m)$ starting at the basepoint 
$e = (e, 0, 0)$ and ending at $(m^n, n, 0)$, defined piecewise by
$$\gamma_n (t) = \left\{ \begin{array}{ll}
					(m^k, k, n(t-k/n)), \,\,\, \frac{k}{n}\leqs t< \frac{k+1}{n},\,\,\, k = 0, \ldots, n-1\\
					\gamma_n (1) = (m^n, n, 0).
				      \end{array}
					\right.
$$  
We call a loop $\alpha\co I \to M_{\infty}$ \emph{normal} (at level $n$) if it is based at $e$ and has the form
$\alpha = \gamma_n \centerdot \widetilde{\alpha} \centerdot \gamma_n^{-1}$, where 
$\widetilde{\alpha} (I) \subset M\cross \{n\} \cross \{0\}$.  Note that if $\alpha$ is normal, then its ``middle third'' $\widetilde{\alpha}$ is uniquely defined.  We will often think of $\widetilde{\alpha}$ as a loop in $M$ rather than in $M_{\infty}(m)$.  

Given a normal loop $\gamma_n \centerdot \widetilde{\alpha} \centerdot \gamma_n^{-1}$, we define its $k$th renormalization to be the normal loop $\gamma_{n+k} \centerdot \widetilde{\alpha} * m^k \centerdot \gamma_{n+k}^{-1}$.
\end{definition}

\begin{lemma}$\label{normalization}$ Let $\gamma\co [0,1]\to M_{\infty} (m)$ be a loop based
at $e$.  Then $\gamma$ is homotopic \emph{(}rel $\{0,1\}$\emph{)} to a normal loop.  
Moreover, every normal loop is homotopic to all of its renormalizations.
\end{lemma} 
\begin{proof}
In general, say we are given a loop $\gamma$ in the mapping cylinder $M_f$ of a map $f\co X\to Y$, and say $\gamma$ is based at $(x_0, 0)\in M_f$ for some $x_0\in X$.  Then if $H_t$ denotes a homotopy from $Id_{M_f}$ to the retraction $r\co M_f \to Y$ (such that $H_t (x_0, 0) = (x_0, t)$ for $t<1$) we have a homotopy connecting $\gamma$ to a loop whose middle third lies in $Y$: 
$$\gamma_t (s) = \left\{ \begin{array}{ll}
					H_{3s} (x_0,0), \,\,\,\,\,\,\,\,\,\,\,\,\,\,\,\,\,\,\,\,\,\,\,\,\,\,\,\,\, 0\leqs s\leqs 														t/3\\
					H_t\left(\gamma\left( \frac{s-t/3}{1-2t/3} \right) \right), \,\,\,
						t/3 \leqs s \leqs 1 - t/3 \\
					H_{3(1-s)} (x_0,0), \,\,\,\,\,\, 1 - t/3 \leqs s \leqs \,1.
				      \end{array}
					\right.
$$  
Now, each loop in $M_{\infty}(m)$ lies in some finite telescope $M_{N} (m)$, and applying the above process $N$ times produces a normal loop (up to reparametrization).  Homotopies between a normal loop and its renormalizations are then produced similarly.
\end{proof}

\begin{definition}  Given loops $\widetilde{\alpha}$ and $\widetilde{\beta}$ in $M$, let $p_{\widetilde{\alpha}, \widetilde{\beta}}$ denote their pointwise sum.
For normal loops $\alpha = \gamma_n\centerdot\widetilde{\alpha}\centerdot\gamma_n^{-1}$ and 
$\beta = \gamma_n\centerdot\widetilde{\beta}\centerdot \gamma_n ^{-1}$ in $M_{\infty} (m)$, 
we define $\alpha*\beta$ to be the normal loop (of level $2n$) given by
$$\alpha*\beta = \gamma_{2n} \centerdot p_{\widetilde{\alpha}, \widetilde{\beta}} \centerdot \gamma_{2n}^{-1}.$$
\end{definition}

\begin{lemma}$\label{properties}$ 
For each pair of normal loops $\alpha$ and $\beta$ of level $n$ (in $M_{\infty} (m)$),  there is basepoint preserving homotopy $\alpha*\beta\heq \beta*\alpha$.
\end{lemma}
\begin{proof} Since $m$ is anchored, there is a homotopy $H\co M\cross M\cross I\to M$ such that
$H(x,y,0) = x* y$, $H(x,y,1) = y* x$, and $H(m^n, m^n, s) = m^{2n}$ for all $s\in I$ and all $n\in \bbN$.  Let 
$\alpha = \gamma_n\centerdot\widetilde{\alpha}\centerdot\gamma_n^{-1}$, 
let $\beta = \gamma_n\centerdot\widetilde{\beta}\centerdot\gamma_n^{-1}$,
and define $h_s (\alpha, \beta)$ to be the loop 
$$h_s (\alpha, \beta) (t) = H (\widetilde{\alpha} (t), \widetilde{\beta}(t), s)$$
(note that $h_s (\alpha, \beta)$ is based at $H(m^n, m^n, s) = m^{2n}$).  The family of loops (based at $e$) given by
$$p_s = \gamma_{2n} \centerdot h_s (\alpha, \beta) \centerdot \gamma_{2n}^{-1}$$   
now provides the desired homotopy between $\alpha*\beta$ and $\beta*\alpha$.
\end{proof}

\begin{proof}[Proof of~\fullref{abelian}]  Let $a$ and $b$ be elements of $\pi_1 (M_{\infty} (m), e)$.  By Lemma~\ref{normalization}, we may choose normal representatives $\alpha = \gamma_n \centerdot  \widetilde{\alpha}\centerdot \gamma_n ^{-1}$ and $\beta = \gamma_k \centerdot \widetilde{\beta}\centerdot \gamma_k^{-1}$ for $a$ and $b$.  Renormalizing if necessary, we can assume $k=n$.  By Lemma~\ref{properties}, it suffices to show that $\alpha\centerdot \beta \heq \alpha * \beta$ (rel \{0,1\}).  Let $m^n$ denote the constant loop at $m^n$.
Note that the $n$th renormalization of $\alpha$
is precisely $\alpha* (\gamma_n \centerdot  m^n \centerdot \gamma_n ^{-1})$ (and similarly for $\beta$).  Using Lemma~\ref{normalization} and commutativity of $*$ we now have
\begin{eqnarray*}
\alpha\centerdot\beta & \heq & (\alpha* (\gamma_n \centerdot  m^n \centerdot \gamma_n ^{-1})) \centerdot (\beta* (\gamma_n \centerdot  m^n \centerdot \gamma_n^{-1})) \\
& \heq & (\alpha* (\gamma_n \centerdot  m^n \centerdot \gamma_n ^{-1})) \centerdot ((\gamma_n \centerdot  m^n \centerdot \gamma_n ^{-1})* \beta)\\
& = &(\gamma_{2n} \centerdot p_{\widetilde{\alpha}, m^n} \centerdot \gamma_{2n}^{-1}) \centerdot (\gamma_{2n} \centerdot p_{m^n, \widetilde{\beta}} \centerdot \gamma_{2n}^{-1})\\  
& \heq &\gamma_{2n} \centerdot p_{\widetilde{\alpha}, m^n} \centerdot p_{m^n, \widetilde{\beta}} \centerdot \gamma_{2n}^{-1} \\
& = &\gamma_{2n}  \centerdot p_{\widetilde{\alpha} \centerdot m^n, m^n \centerdot \widetilde{\beta}}  \centerdot \gamma_{2n}^{-1}\\
& = &(\gamma_n \centerdot  (\widetilde{\alpha} \centerdot m^n)\centerdot  \gamma_{n}^{-1})* (\gamma_{n} \centerdot (m^n \centerdot \widetilde{\beta}) \centerdot \gamma_{n}^{-1}).
\end{eqnarray*}
Since $\widetilde{\alpha} \centerdot m^n \heq \widetilde{\alpha}$,
we have a homotopy $\gamma_n \centerdot \widetilde{\alpha}_s \centerdot \gamma_n^{-1}$ from $\gamma_{n}\centerdot (\widetilde{\alpha} \centerdot m^n) \centerdot\gamma_{n}^{-1}$ to $\alpha$ (we may assume each loop in this homotopy is normal) and analogously for $\beta$.
The family of loops $\gamma_{n}\centerdot \widetilde{\alpha}_s \centerdot\gamma_{n}^{-1} * 
\gamma_{n}\centerdot \widetilde{\beta}_s \centerdot\gamma_{n}^{-1}$
provides a homotopy from $(\gamma_n \centerdot  (\widetilde{\alpha} \centerdot m^n)\centerdot  \gamma_{n}^{-1})* (\gamma_{n} \centerdot (m^n \centerdot \widetilde{\beta}) \centerdot \gamma_{n}^{-1})$ 
to $\alpha * \beta$,
and since all homotopies involved are basepoint preserving, this completes the proof.
\end{proof}

We now show that all components of $M_{\infty} (m)$ have abelian fundamental group, not just the component containing $e$.

\begin{corollary}$\label{abelian2}$  Let $M$ be a homotopy commutative monoid which is stably group-like with respect to an anchored element $m\in M$.  Then all path components of $M_{\infty} (m)$ have abelian fundamental group.
\end{corollary}
\begin{proof}  For each element $(x, n, t)\in M_{\infty} (m)$, we define $C_{(x,n,t)}$ to be the component of $M_{\infty} (M)$ containing this element.  We need to show that
$\pi_1 (C_{(x,n,t)} (M))$ is abelian.  Let $x^{-1}\in M$ be a stable homotopy inverse for $x$, i.e. an element such that for some $N$, $x* x^{-1}$ and $m^N$ lie in the same component of $M$.  Note that by adding $m^n$ to $x^{-1}$ if necessary, we may assume that $N>n$.  We will construct maps 
$$f\co (C_{(x,n,t)}, (x,n,t))\to (C_{(e,0,0)}, (x^{-1}* x, N, t))$$
and
$$g\co (C_{(e,0,0)}, (x^{-1}* x, N, t)) \to (C_{(x,n,t)}, (x* x^{-1}* x, n+N, t))$$
and show that the composition $g_*\circ f_*$ is injective on $\pi_1$, from which it follows that $f_*$ is injective.  This will suffice, since by Proposition~\ref{abelian} the group $\pi_1 (C_{(e,0,0)}, (x^{-1}* x, N, t))$ is abelian.

The maps $f$ and $g$ are defined by 
\begin{eqnarray*}
f(y, k, s) &= &(x^{-1} * y, k+(N-n), s)\\
g(y,k,s) &=& (x* y, k+n, s);
\end{eqnarray*}
note that in both cases these are continuous maps (defined, in fact, on the whole mapping telescope $M_{\infty} (m)$) and they map the basepoints in the manner indicated above.  The composite map
is given by $g\circ f (y,k,t) = (x* x^{-1}* y, k + N, t)$.

Consider an element $[\alpha]\in \ker(g_*\circ f_*)$.  Then $\alpha$ lies in some finite telescope
$M_k (m)$, and by collapsing this telescope to its final stage, we obtain a \emph{free} homotopy from 
$\alpha$ to a loop $\overline{\alpha}$ lying in $M\cross \{k\} \cross \{0\}$.  Now, we have a free homotopy $$g\circ f\circ \alpha \heq g\circ f\circ \overline{\alpha} = x* x^{-1} * \overline{\alpha}$$
where the final loop lies in $M\cross \{k+N\}\cross \{0\}$.  By assumption, there is a path in $M$ from $x* x^{-1}$ to $m^N$, and together with homotopy commutativity of $M$ we find that
$x* x^{-1} * \overline{\alpha}\heq m^N* \overline{\alpha} \heq \overline{\alpha}* m^N.$
But this loop, lying in $M\cross \{k+N\}\cross \{0\}$, is homotopic (in $M_{k+N} (m)$) to the loop $\overline{\alpha}$ lying in $M\cross \{k\}\cross \{0\}$ (the homotopy is obtained by collapsing
the telescope $M_{k+N} (m)$ to its final stage).  
By construction, $\overline{\alpha}$ is homotopic to $\alpha$.  Thus we have a free homotopy
$g\circ f\circ \alpha \heq \alpha$, and by assumption $g\circ f\circ \alpha$ is nullhomotopic.  
Hence $\alpha$ is \emph{freely} nullhomotopic.  But
freely nullhomotopic loops are always trivial in $\pi_1$, so $g_*\circ f_*$ is injective as claimed.
\end{proof}

\begin{proof}[Proof of~\fullref{gp-completion-thm}] 
Recall that by the McDuff--Segal Theorem, we have a natural zig-zag
\begin{equation}\label{gc-zig-zag2}
M_{\infty} (m) \stackrel{i_m}{\maps} \hofib (q_m) \stackrel{\heq}{\longleftarrow} \Omega BM,
\end{equation}
where the first map induces an isomorphism in homology for every local (abelian) coefficient system,
and the second map is a weak equivalence.
By Corollary~\ref{abelian2}, all components of $M_{\infty} (m)$ have abelian fundamental group, and hence $i_m$ induces isomorphisms on $\pi_1 \isom H_1$ (note that $\pi_1 \hofib (q_m) \isom \pi_1 \Omega BM$ is abelian since $\Omega BM$ is an $H$--space).  It is well-known that a map inducing isomorphisms on homology with local coefficients, and on $\pi_1$, is a weak equivalence (see, for example, Hatcher~\cite[p. 389, Ex. 12]{Hatcher}).  

To complete the proof of Theorem~\ref{gp-completion-thm}, we must show that the zig-zag (\ref{gc-zig-zag2}) induces an isomorphism of groups $\pi_0 (M_{\infty} (m)) \isom \pi_0(\Omega BM)$ (the multiplication on $\pi_0  (M_{\infty} (m))$ was described in Proposition~\ref{stbly-gp-like}).  We already know that these maps induce a bijection, so it suffices to check that the induced map is a homomorphism.  
Every component of $M_{\infty} (m)$ is represented by a point of the form $(x, n, 0)$, with $x\in M$ and
$n\in \bbN$.  Now, the fiber of $q_m$ over $*\in BM$ is precisely the objects of the category $(M_{\infty})_M$, i.e. the space $M_{\infty} (m)$, and hence we identify $(x, n, 0)$ with a point in
$q_m^{-1} (*)$.  Hence we may write $i_m (x, n, 0) = \left( (x, n, 0), c_* \right) \in \hofib (q_m)$, where
$c_*$ denotes the constant path at $*\in BM$.
Next, recall that since $BM$ is the realization of a category with $M$ as morphisms, every element $y\in M$ determines a loop $\alpha_y \in \Omega BM$.  Let $\psi$ denote the natural map from $\Omega BM \to \hofib (q_m)$.  We claim the points
$\psi \left( \alpha_{m^n}^{-1} \centerdot \alpha_x \right)$ 
and 
$\left( (x, n, 0), c_* \right)$ lie in the same path component of $\hofib (q_m)$.  This implies that 
the map $\pi_0 (M_{\infty} (m)) \to \pi_0(\Omega BM)$ sends the component of $(x,n,0)$ to the component of $\alpha_{m^n}^{-1} \centerdot \alpha_x$.   Since $M$ is homotopy commutative
and $\pi_0 (\Omega BM)$ is the group completion of $\pi_0 M$, this map is a homomorphism of monoids.

We now produce the required path (in $\hofib (q_m)$) between the points $\left( (x, n, 0), c_* \right)$ and
$\psi (\alpha_{m^n}^{-1} \centerdot \alpha_x )$.  By definition of the map $\psi$ we have
$\psi ( \alpha_{m^n}^{-1} \centerdot \alpha_x ) = (e, \alpha_{m^n}^{-1} \centerdot \alpha_x )$, where $e = (e, 0, 0) \in M_{\infty} (m)$.  There are morphisms in the category 
$(M_{\infty})_M$ from the object $(e, n, 0)$ to $(x, n, 0)$ and to $(m^n, n, 0)$, corresponding (respectively) to the elements $x$ and $m^n$ in $M$.
These morphisms give paths $\beta_x$ and $\beta_{m^n}$ in $\abs{(M_{\infty})_M}$ which map under
$q_m$ to the paths $\alpha_x$ and $\alpha_{m^n}$, respectively.  Letting $\alpha_x^t$ denote the path
$\alpha_x^t (s) = \alpha_x (1-t + ts)$, the formula $t\mapsto (\beta_x (1-t), \alpha_x^t)$ defines a path in $\hofib (q_m)$ starting at $\left( (x, n, 0), c_*\right)$ and ending at $\left( (e, n, 0), \alpha_x\right)$.  One next constructs an analogous path from $\left( (e, n, 0), \alpha_x \right)$ to 
$\left( (m^n, n, 0), \alpha_{m^n}^{-1} \centerdot \alpha_x \right)$.  Finally, since $(m^n, n, 0)$ and $e$ lie in the same component of $M_{\infty} (m) = q_m^{-1} (*)$, we have a path in 
$\hofib(q_m)$ from $\left( (m^n, n, 0), \alpha_{m^n}^{-1} \cdot \alpha_x \right)$ to $(e, \alpha_{m^n}^{-1} \cdot \alpha_x)$.
\end{proof}

\section{Group completion in deformation $K$--theory}$\label{gp-completion-2}$

We now apply the results of the previous section to deformation $K$--theory. We work mainly in the unitary case, but all of the results are valid in the general linear case as well (and we have noted the places in which the arguments differ).

For applications to deformation $K$--theory, our real interests lie in the monoid of homotopy orbit spaces
$$\Rep(G)_{h U} = \coprod_{n=0}^{\infty} \Hom(G, U(n))_{hU(n)},$$
but we can often work with the simpler monoid $\Rep(G)$ instead.
Note that the spaces $EU(n)$ are naturally based: recall from Section~\ref{K-def} that $EU(n)$ is the classifying space of a topological category whose object space is $U(n)$; the object $I_n\in U(n)$ provides the desired basepoint $*_n\in EU(n)$.  These basepoints behave correctly with respect to block sums, i.e. $*_n\oplus *_m = *_{n+m}$.

\begin{lemma}$\label{reduction1}$
Let $G$ be a discrete group.  Then $\Rep(G)$ is stably group-like with respect to $\psi\in \Hom(G, U(m))$ if and only if $\Rep(G)_{hU}$ is stably group-like with respect to $[*_m, \psi]\in  \Hom(G, U(m))_{hU(m)}$.
\end{lemma}
\begin{proof}  Say $\Rep(G)$ is stably group-like with respect to $\psi\in \Hom(G, U(m))$.  Then given a point $[e, \rho]\in \Rep(G)_{hU}$ (with $e\in EU(n)$ and $\rho\co G\to U(n)$ for some $n$), we know that there is a representation $\rho^{-1}\co G\to U(k)$ such that $\rho\oplus \rho^{-1}$ lies in the component of the $\psi^{l}$ (the $l$-fold block sum of $\psi$ with itself), where $l = \frac{n+k}{m}$.  Now, for each $e'\in EU(k)$, the point $[e', \rho^{-1}]$ is a stable homotopy inverse for $[e, \rho]$ (with respect to $\psi$), since there is a path in $EU(n+k)\cross \Hom(G, U(n+k))$ from $[e\oplus e', \rho\oplus \rho^{-1}]$ to $[*_{n+k}, \psi^l]$.

Conversely, if $\Rep(G)_{hU}$ is stably group-like with respect to $[*_m, \psi]$, then each
element $[e, \rho]\in \Hom(G, U(n))_{hU(n)}$
has a stable homotopy inverse $$[e', \rho^{-1}]\in \Hom(G, U(k))_{hU(k)}$$ (for some $k$),
i.e. there is a path in $\Hom(G, U(n+k))_{hU(n+k)}$ from 
$[e\oplus e', \rho\oplus \rho^{-1}]$ to $[*_{n+k}, \psi^l]$ (where again $l = \frac{n+k}{m}$).  Path-lifting for the fibration $EU(n+k)\cross \Hom(G, U(n+k))\to \Hom(G, U(n+k))_{hU(n+k)}$
produces a path in $EU(n+k)\cross \Hom(G, U(n+k))$ from $(e\oplus e', \rho\oplus \rho^{-1})$ to some point
$(*_{n+k}\cdot A, A^{-1} \psi^l A)$, with $A\in U(n+k)$.  The second coordinate of this path, together with connectivity of $U(n)$, shows that $\rho\oplus\rho^{-1}$ lies in the component of $\psi^l$, i.e. $\rho^{-1}$ is in fact a stable homotopy inverse for $\rho$ (with respect to $\psi$).
\end{proof}

We will now show that in the monoid $\Rep(G)_{hU}$, elements are always anchored.  First we need some lemmas regarding the unitary and general linear groups, which are probably well-known.

\begin{lemma}$\label{centralizer1}$
Consider an element $D = \lambda_1 I_{n_1} \oplus \cdots \oplus \lambda_k I_{n_k}\in GL_n (\bbC)$, where $n=\sum n_i$ and the $\lambda_i$ are distinct.  Then the centralizer of $D$ in $GL_n (\bbC)$ is the subgroup 
$GL(n_1)\cross\cdots \cross GL(n_k),$
embedded in the natural manner.  As a consequence, the analogous statement holds for the unitary groups.
\end{lemma}

\begin{lemma}$\label{centralizer2}$
Let $K\subset GL_n (\bbC)$ be a subgroup.  Then the set of diagonalizable matrices in the centralizer $C(K)$ is connected.  Similarly, for each $K\subset U(n)$, the centralizer of $K$ in $U(n)$ is connected.
\end{lemma}
\begin{proof}  We prove the general linear case; the argument for $U(n)$ is nearly identical (since by the Spectral Theorem every element of $U(n)$ is diagonalizable).  Let $A\in C(K)$ be diagonalizable.  We will produce a path (of diagonalizable matrices) in $C(K)$ from $A$ to the identity.  Choose $X\in GL_n (\bbC)$ such that
$XAX^{-1} =  \lambda_1 I_{n_1} \oplus \cdots \oplus \lambda_k I_{n_k}$
(for some $n_i$ with $\sum n_i = n$).  Then by Lemma~\ref{centralizer1} we have
$XKX^{-1}\subset GL(n_1)\cross\cdots \cross GL(n_k)$.
Now, choose paths $\lambda_i (t)$ from $\lambda_i$ to $1$, lying in $\bbC-\{0\}$ (or in the unitary case, lying in $S^1$).  This gives a path of matrices $Y_t$ connecting $XAX^{-1}$ to $I$, and for each $t\in I$ we have $Y_t \in C(XKX^{-1})$.  Now $X^{-1} Y_t X$ is a path from $A$ to $I$ lying in $C(K)$.
\end{proof}

\begin{corollary}$\label{gp-completion-cor}$
Let $G$ be a finitely generated discrete group such that $\Rep(G)$ is stably group-like with respect to a representation $\rho \in \Hom(G,U(k))$.
Then there is a natural isomorphism 
$$\pi_* \K(G) \homeo \pi_* \tele \left( \Rep(G, U)_{hU} \stackrel{\oplus \rho}{\maps} \Rep(G, U)_{hU} 
			\stackrel{\oplus\rho}{\maps} \cdots \right),
$$
where $\oplus \rho$ denotes block sum with the point $[*_k, \rho]\in \Hom(G, U(k)_{hU(k)}$.  The analogous statement holds for general linear deformation $K$--theory.
\end{corollary}

\begin{remark}   Naturality here has essentially the same meaning as in Theorem~\ref{gp-completion-thm}, i.e. these spaces are connected by a zig-zag of natural weak equivalences.  The comments after Theorem~\ref{gp-completion-thm} regarding basepoints apply here as well.

Various examples of groups to which Corollary~\ref{gp-completion-cor} applies directly are discussed in the author's thesis~\cite[Chapter 6]{Ramras-thesis}.  These groups include the fundamental group of any surface, and groups admitting similar presentations, as well as finitely generated abelian groups. 
\end{remark}

\begin{proof}[Proof of~\fullref{gp-completion-cor}] The result will follow immediately from Lemma~\ref{zeroth-space}, Proposition~\ref{classifying-space},  Lemma~\ref{reduction1}, and Theorem~\ref{gp-completion-thm} once we show that the element 
$$[*_k, \rho]\in \Hom(G, U(k))_{hU(k)}$$
is anchored in the monoid $\Rep(G)_{hU} = \abs{\mcR(G)}$.  We will work with $\abs{\mcR(G)}$; note that the element $[*_k, \rho]$ above maps to the object 
$\rho\in \Hom(G, U(k))$ under the isomorphism $\alpha$ 
of monoids in Lemma~\ref{classifying-space}.

Given a collection of matrices $X = \{X_{(n,m)}\}_{n,m\in \bbN}$ with $X_{(n,m)}\in U(n+m)$, we can define a functor $F_X\co \mcR(G)\cross \mcR(G) \to \mcR(G)$
as follows. 
Given objects $\psi_1\in \Hom(G, U(n))$ and $\psi_2\in \Hom(G, U(m))$, we set
$$F_X (\psi_1, \psi_2) = X_{(n,m)} (\psi_1\oplus \psi_2) X_{(n,m)}^{-1}.$$
We define $F_X$ on morphisms by sending
$(A, B)\co (\psi_1, \psi_2)\to (A \psi_1 A^{-1}, B \psi_2 B^{-1})$
to the morphism 
$$X_{(n,m)} (\psi_1 \oplus \psi_2) X_{(n,m)}^{-1}
		\maps X_{(n,m)} \left( (A \psi_1 A^{-1}) \oplus (B \psi_2 B^{-1}) \right) X_{(n,m)}^{-1}$$
represented by the matrix $X_{(n,m)} (A\oplus B) X_{(n,m)}^{-1}$.

Let $\tau_{n,m}$ be the matrix
\vspace{.02in}
$$\left[ 
\begin{array}{cc}
       0 & I_m \\           
       I_n & 0 
\end{array}
\right]\vspace{.03in}
$$
and choose paths $\gamma_{n,m}$ from $I_{n+m}$ to $\tau_{n,m}$ in $U(n+m)$.  When $n = m=kl$ ($l\in \bbN$) we may assume, by Lemma~\ref{centralizer2}, that $\gamma_{kl,kl} (t)\in \mathrm{Stab} (\rho^{2l})$ for all $t\in I$ (note that $\Stab(\rho^{2l}) = C (\mathrm{Im}\rho^{2l})$).
Let $X^t$ denote collection $X^t_{n,m} = \gamma_{n,m} (t)$, and let $F_t = F_{X^t}$ be the associated functor.  Then $F_0 = \oplus$ is the functor inducing the monoid structure on $\abs{\mcR(G)}$, i.e. $\abs{F_0} (x, y) = x\oplus y$, and $\abs{F_1} (x, y) = y\oplus x$.  Moreover, at every time $t$ we have
$\abs{F_t}(\rho^l, \rho^l) = \rho^{2l}$.  The path of functors $F_t$ provides the desired homotopy, proving that $\rho\in \abs{\mcR(G)}$ is anchored.  (Note that a continuous family of functors $G_t \co \mcC\to \mcD$ defines a continuous functor $\mcC \cross \mathcal{I}\to \mcD$, where $\mathcal{I}$ denotes the topological category whose object space and morphism space are both the unit interval $[0,1]$, and hence yields a continuous homotopy).
\end{proof}

\begin{remark}
We note that in the above proof there are natural isomorphisms between $F_0$ and $F_1$, given by the matrices $\tau_{n,m}$.  This is the usual way to show that a monoid coming from a permutative category is homotopy commutative, but this homotopy does not anchor $\rho$.
\end{remark}


\section{Excision for free products}$\label{free-excision}$

In this section we present our results on the excision problem for free products and discuss some resulting computations.  First we describe the excision problem more generally, in the context of amalgamated products.

Let $G$, $H$, and $K$ be discrete groups, with homomorphisms $f_1\co K\to G$ and 
$f_2\co K\to H$.  Then associated to the co-cartesian (i.e. pushout) diagram of groups
$$\xymatrix{
       {K} \ar[r]^{f_1} \ar[d]^{f_2}
       & {G} \ar[d] \\           
       {H} \ar[r] 
       & {G *_K H}  
                     }
$$
there is a diagram of spectra
\begin{equation}\label{excision-diagram}
\xymatrix{
       {\K (G *_K H)} \ar[r] \ar[d]
       & {\K (G)} \ar[d]^{f_1^*} \\           
       {\K (H)} \ar[r]^{f_2^*} 
       & {\K (K).}  
                     }
\end{equation}
We will say that the amalgamated product $G*_K H$ satisfies excision (for deformation $K$--theory) if diagram (\ref{excision-diagram})
is homotopy cartesian, i.e. if the natural map from $\K (G *_K H)$ to the homotopy pullback is a weak equivalence.  Note that since we are dealing with connective $\Omega$--spectra, this is the same as saying that the diagram of zeroth spaces is homotopy cartesian.

Excision results are important from the point of view of computations, since associated to a homotopy cartesian diagram of spaces
$$\xymatrix{
       {W} \ar[r]^{f} \ar[d]^{g}
       & {X} \ar[d]^{h} \\           
       {Y} \ar[r]^{k} 
       & {Z}  
                     }
$$
there is a long exact ``Mayer--Vietoris'' sequence of homotopy groups
\begin{equation}\label{LES}
 \ldots \maps \pi_k (W)     \stackrel{f_* \oplus g_*}{\maps} 
     \pi_k (X)\oplus \pi_k (Y)  \stackrel{h_* - k_*}{\maps} 
     \pi_k (Z) \stackrel{\partial}{\maps} 
     \pi_{k-1} (W)\maps \ldots
\end{equation}
which comes from combining the long exact sequences associated to the vertical maps (see Hatcher~\cite[p. 159]{Hatcher}; note that the homotopy fibers of the vertical maps in a homotopy cartesian square are weakly equivalent).  It is not difficult to check that if all the spaces involved are group-like $H$--spaces, and the maps are homomorphisms of $H$--spaces, then the maps in this sequence (including the boundary maps) are homomorphisms in dimension zero.  Hence when applied to (the zeroth spaces of) the deformation $K$--theory in an amalgamation diagram, assuming excision one obtains a long exact sequence in $\K^*$.  

Deformation $K$--theory can fail to satisfy excision in low dimensions.  As mentioned in the introduction, such examples arise from connected sum decompositions of Riemann surfaces (see Ramras~\cite[Section 6]{Ramras-surface}).
Using the results from Section~\ref{gp-completion-2}, we will show in Theorem~\ref{excision-free-prod} that deformation $K$--theory satisfies excision for free products.

The proof of Theorem~\ref{excision-free-prod} requires several lemmas.  

\begin{lemma}$\label{cartesian-orbit}$
For all discrete groups $G$, $H$ and $K$, the homotopy orbit space $\Hom(G*_K H, U(n))_{hU(n)}$
is naturally homeomorphic to the pullback
$$ \Hom(G, U(n))_{hU(n)} \cross_{\Hom(K, U(n))_{hU(n)}} \Hom(H,U(n))_{hU(n)}.$$
The analogous statement holds for the general linear groups in place of the unitary groups.
\end{lemma}
\begin{proof}  It follows from the proof of Proposition~\ref{classifying-space} that for all groups $L$, the space $\Hom(L, U(n))_{hU(n)}$ is homeomorphic to the realization of a simplicial space of the form
$$\abs{k\mapsto U(n)^k \cross \Hom(L, U(n))}.$$
The lemma now follows from the fact that geometric realization commutes with pullbacks in the category of compactly generated spaces (see Gabriel--Zisman~\cite[III.3]{Gabriel-Zisman}).  The proof for $GL_n (\bbC)$ is identical.
\end{proof}

We will need to consider certain submonoids of the topological monoid $\Rep(G)_{hU}$ underlying deformation $K$--theory.
Given a topological monoid $M$ and a submonoid $A\subset \pi_0(M)$, we may define a corresponding submonoid $M(A)\subset M$ by setting
$$M(A) = \bigcup_{C\in A} C.$$
For every such $A$, we have
$\pi_0 M(A) = A$, and $M(A)$ is a union of connected components of $M$.

\begin{lemma}$\label{M(A)}$
Let $(M, *)$ be a homotopy commutative topological monoid, and let $m\in M$ be anchored.
Consider a submonoid $A\subset \pi_0 (M)$ containing the component of $m$.
Then $m$ is anchored in the submonoid $M(A) \subset M$.
\end{lemma}
\begin{proof} If $H\co M\cross M\cross I \to M$ is a homotopy anchoring $m$, then we claim that
$H(M(A)\cross M(A) \cross I) \subset M(A)$.  Indeed, if $x$ and $y$ lie in $M(A)$, then for every
$t\in I$, the points $H(x,y,t)$ and $H(x, y, 0)$ lie in the same component of $M$, and $H(x,y,0) = x*y$ lies in $M(A)$ because $A$ is a submonoid of $\pi_0 M$.
\end{proof}

In order to apply the results of Section~\ref{gp-completion-2}, we will need to filter the monoid $\Rep(G)_{hU}$ by submonoids which are stably group-like with respect to some representation.  We will need a lemma regarding the algebraic nature of representation spaces.  This lemma requires the group $G$ to be finitely generated, and is the reason for this assumption in Theorem~\ref{excision-free-prod}.

Recall from Whitney~\cite{Whitney} that a real algebraic variety $V\subset \R^n$ is the set of common zeros of some collection of polynomials with real coefficients.

\begin{lemma}$\label{real-alg}$
Let $G$ be a finitely generated discrete group.  Then the representation spaces
$\Hom(G, U(n))$ and $\Hom(G, GL_n(\bbC))$ are homeomorphic to real algebraic varieties.
\end{lemma}
\begin{proof} First, note that $U(n)$ is a real algebraic variety, since
we may write
$$U(n) = \{A\in M_n (\bbC)\isom \R^{2n^2} \,:\, AA^* = I_n\},$$
and the equation $AA^* = I_n$ is simply a system of polynomial equations.
Similarly, we have
$$GL_n (\bbC) \isom \{(A,B)\in M_n(\bbC)^2 \homeo \R^{4n^2} \,:\, AB = I\}.$$
Now if $G = \langle g_1, \ldots, g_k | \{r_i\}_{i\in \mcR} \rangle$ 
is a presentation for $G$, then we may write
\begin{equation*}\begin{split} 
\Hom(G, U(n)) \isom \{ (A_1, \ldots, A_k)\in M_n (\bbC)^k \,: \, &A_j A_j^* = I_n \,\,(j= 1, \ldots, k) \\
&\mathrm{and} \,\,r_i(A_1, \ldots, A_k) = I_n \,\,(i\in \mcR) \}.
\end{split}\end{equation*}			     
			    		
The fact that this topology on $\Hom(G, U(n))$ agrees with the compact-open topology follows from the fact that equality of two representations on a generating set implies equality, so convergent sequences of representations are simply those converging on the chosen generators.  The situation for the general linear representation varieties is similar.
\end{proof}								

\begin{corollary}$\label{finite-comp}$
If $G$ is a finitely generated discrete group, then the representation spaces
$\Hom(G, U(n))$ and $\Hom(G, GL_n (\bbC))$ have finitely many connected components.
\end{corollary}
\begin{proof}  By Lemma~\ref{real-alg}, these spaces are real algebraic varieties, hence triangulable (see Hironaka~\cite{Hironaka}).  Hence their connected components and their path components coincide.  The unitary representation varieties are compact, since they are closed subsets of a product of copies of $U(n)$.   Hence they have finitely many connected components, and also finitely many path components.  More generally, Whitney's theorem~\cite{Whitney} states that every real algebraic variety has finitely components, so the result follows in the general linear case as well.
\end{proof}

We can now prove the main result of this section.

\begin{theorem}$\label{excision-free-prod}$
Let $G$ and $H$ be finitely generated discrete groups.
Then the diagram of spectra 
$$\xymatrix{
       {\K(G * H)} \ar[r] \ar[d]
       & {\K(G)} \ar[d]\\           
       {\K(H)} \ar[r] 
       & {\K(\{1\})}  
                     }
$$ 
is homotopy cartesian.
\end{theorem}

Note that $\K (\{1\})\heq \mathbf{ku}$, the complex connective $K$--theory spectrum, since (in the unitary case) 
$\Rep(\{1\})_{hU} = \coprod_{n=0}^{\infty} BU(n)$.  This also holds for general linear deformation $K$--theory, because $U(n)\heq GL_n (\bbC)$.

\begin{remark}$\label{cart-htpy-gps}$
Theorem~\ref{excision-free-prod} yields a long exact sequence in $\K^*$, and the boundary maps in this sequence are always zero because the map $\K^* (G) \to \K^* (\{1\})$ induced by the inclusion
$\{1\}\injects G$ is split by the map $\K^*(\{1\}) \to \K^*(G)$ induced by the projection $G\to \{1\}$.
Hence Theorem~\ref{excision-free-prod} produces cartesian diagrams of homotopy groups in each dimension, and in odd dimensions we have $\K^* (\{1\}) = \pi_* \mathbf{ku} = 0$, meaning that $\K^* (G*H) = \K^* (G) \oplus \K^* (H)$.
\end{remark}

\begin{proof}[Proof of~\fullref{excision-free-prod}]
The proofs for the general linear and unitary cases are identical, so we work in the unitary case.  The proof involves reducing to a diagram of homotopy orbit spaces, which will be homotopy cartesian by Lemma~\ref{cartesian-orbit}.

We begin by noting that since the spectra involved are all connective $\Omega$-spectra, it will suffice to show that the diagram of zeroth spaces is homotopy cartesian.  In order to apply Theorem~\ref{gp-completion-thm} we need to filter the underlying monoids by submonoids which are stably group-like with respect to compatible representations. 
For each $n=1, 2, \ldots$, we will define submonoids $\Rep(G*H)_{hU}^{(n)}\subset \Rep(G*H)_{hU}$, $\Rep(G)_{hU}^{(n)}\subset \Rep(G)_{hU}$ and $\Rep(H)_{hU}^{(n)}\subset \Rep(H)_{hU}$ having the following properties:
\begin{enumerate}
\item Each of these submonoids is of the form $\Rep(\cdot)_{hU} (A)$ for
some submonoid $A\subset \pi_0 \Rep(G)_{hU}$ (see Lemma~\ref{M(A)}), and
 $\Rep(\cdot)_{hU} = \bigcup_n \Rep(\cdot)_{hU}^{(n)}$.

\item Under the natural maps from $\Rep(G*H)_{hU}$ to $\Rep(G)_{hU}$ and $\Rep(H)_{hU}$, $\Rep(G*H)_{hU}^{(n)}$ maps to $\Rep(G)_{hU}^{(n)}$ and to $\Rep(H)_{hU}^{(n)}$ respectively.

\item There are representations $\rho_n$ of $G$ and $\psi_n$ of $H$ (of the same dimension $d=d(n)$) such that $\Rep(G*H)_{hU}^{(n)}$ is stably group-like with respect to $[*_d, (\rho_n, \psi_n)]$ and $\Rep(G)_{hU}^{(n)}$ and $\Rep(H)_{hU}^{(n)}$ are stably group-like with respect to $[*_d, \rho_n]$ and $[*_d, \psi_n]$ (respectively).

\item For each $n$, the square
$$\xymatrix{
       {\Rep(G * H)_{hU}^{(n)}} \ar[r] \ar[d]
       & {\Rep(G)_{hU}^{(n)}}  \ar[d] \\           
       {\Rep(H)_{hU}^{(n)}} \ar[r]
       & {\Rep(\{1\})_{hU}}  
                     }
$$
is cartesian, i.e. the natural map 
$$\Rep(G * H)_{hU}^{(n)} 
	\maps \lim \left( \Rep(G)_{hU}^{(n)} \to \Rep(\{1\})_{hU} \leftarrow \Rep(H)_{hU}^{(n)} \right)
$$
is a homeomorphism.
\end{enumerate}

Assuming the existence of such filtrations, we now complete the proof of Theorem~\ref{excision-free-prod}.  By Lemma~\ref{zeroth-space}, it suffices to show that the diagram
\begin{equation}\label{loop-lim}
\xymatrix{
       {\Omega B(\Rep(G*H)_{hU})} \ar[r] \ar[d]
       & {\Omega B(\Rep(G)_{hU})}  \ar[d] \\           
       {\Omega B(\Rep(H)_{hU})} \ar[r]
       & { \Omega B(\Rep(\{1\})_{hU})}  
                     }
\end{equation}
is homotopy cartesian.
By Property (2) of the filtrations, the diagram
\begin{equation}\label{loop}
\xymatrix{
       {\Omega B(\Rep(G*H)_{hU}^{(n)})} \ar[r] \ar[d]
       & {\Omega B(\Rep(G)_{hU}^{(n)})}  \ar[d] \\           
       {\Omega B(\Rep(H)_{hU}^{(n)})} \ar[r]
       & { \Omega B(\Rep(\{1\})_{hU})}  
                     }
\end{equation}
is well-defined for every $n$.
By Property (1) of the filtrations, we have $\Rep(\cdot)_{hU} = \colim_{n} \Rep(\cdot)^{(n)}_{hU}$.
Now, direct unions commute with geometric realization because both are quotient constructions, so
the classifying space functor $B$ commutes with direct unions of topological monoids.
Also, if $K$ is compact and $X_1\to X_2 \to \cdots$ is a sequence of Hausdorff spaces, then each map from $K$ into $\colim_n X_n$ lands in one of the $X_n$, so we have a homeomorphism
$$\colim_{n\to \infty} \Map (K, X_n) \stackrel{\isom}{\maps} \Map_* (K, \colim X_n).$$
These two facts imply that
Diagram (\ref{loop-lim}) is the colimit (as $n$ tends to infinity) of the diagrams (\ref{loop}).
Hence it will suffice to show that Diagram (\ref{loop}) is homotopy cartesian for each $n$.  

Now, by Property (3) we know that there are representations $\rho_n\co G\to U(d(n))$ and $\psi_n\co H\to U(d(n))$ such that these monoids are stably group-like with respect to the points $[*_{d(n)}, (\rho_n, \psi_n)]$, $[*_{d(n)}, \rho_n]$ and $[*_{d(n)}, \psi_n]$ (respectively).
Furthermore, Lemma~\ref{M(A)}, together wth Property (1), shows that these basepoints are anchored.  So we may apply Theorem~\ref{gp-completion-thm}.

To simplify notation, we let $X^{(n)} = \coprod_{k=0}^{\infty} BU(k)$, $Y^{(n)} = \Rep(G)_{hU}^{(n)}$,
$Z^{(n)} = \Rep(H)_{hU}^{(n)}$, and $W^{(n)} = \Rep(G*H)_{hU}^{(n)}$.  Also, let 
$$W_{\infty}^{(n)} = \colim \left( W^{(n)}\xrightarrow{\oplus [*_{d(n)}, (\rho_n, \psi_n)]} W^{(n)} 
					\xrightarrow{\oplus [*_{d(n)}, (\rho_n, \psi_n)]} \cdots \right)
$$ 
and let
$\widetilde{W}_{\infty}^{(n)}$ denote the homotopy colimit of the same sequence.  We define $X_{\infty}^{(n)}$, $\widetilde{X}_{\infty}^{(n)}$, $Y_{\infty}^{(n)}$, $\widetilde{Y}_{\infty}^{(n)}$, $Z_{\infty}^{(n)}$, and $\widetilde{Z}_{\infty}^{(n)}$ analogously; the direct system for $X$ uses block sum with $*_{d(n)}\in BU(d(n))$.

With this notation, Corollary~\ref{gp-completion-cor} provides a commutative zig-zag of weak equivalences
linking Diagram (\ref{loop}) to the diagram
\begin{equation}\label{tele}
\xymatrix{
       {\widetilde{W}_{\infty}^{(n)}} \ar[r] \ar[d]
       & {\widetilde{Y}_{\infty}^{(n)}}  \ar[d] \\           
       {\widetilde{Z}_{\infty}^{(n)}} \ar[r]
       & {\widetilde{X}_{\infty}^{(n)}.} 
                 }
\end{equation}
Hence for each $n$, Diagram (\ref{loop}) is homotopy cartesian if and only if Diagram (\ref{tele}) is homotopy cartesian.  
The fact that Diagram (\ref{tele}) is homotopy cartesian essentially follows from Property (4) of the filtrations together with the general fact that homotopy pull-backs commute with directed homotopy colimits.  In this case, though, we can provide the following direct argument.

We must show that the natural map
$$\widetilde{W}_{\infty}^{(n)}
\maps \holim \left(\widetilde{Y}^{(n)}_{\infty} \maps \widetilde{X}^{(n)}_{\infty}
						\longleftarrow \widetilde{Z}^{(n)}_{\infty} \right)
$$
is a weak equivalence.  But this map fits into the commutative diagram
\begin{equation}\label{tele-holim}
\xymatrix{
       {\widetilde{W}_{\infty}^{(n)}}\ar[r] \ar[d]^(.44){\begin{turn}{270} $\heq$ \end{turn} }
       & {\holim \left( \widetilde{Y}^{(n)}_{\infty} \maps \widetilde{X}^{(n)}_{\infty}
			              \longleftarrow \widetilde{Z}^{(n)}_{\infty} \right)}  \ar[dd]^(.48){\begin{turn}{270} 																	$\heq$ \end{turn} } \\           
       {W_{\infty}^{(n)}} \ar[d]^(.4){\begin{turn}{270} $\homeo$ \end{turn} }\\
       {\lim \left(Y^{(n)}_{\infty} \maps X^{(n)}_{\infty}
						\longleftarrow Z^{(n)}_{\infty} \right)} \ar[r]^(.47){\alpha}
       & {\holim \left( Y^{(n)}_{\infty} \maps X^{(n)}_{\infty}
						\longleftarrow Z^{(n)}_{\infty} \right).}
                     }
\end{equation}
The maps labeled $\heq$ are weak equivalences because they arise from collapsing mapping telescopes.  Property (4) states that $W^{(n)} \homeo \lim (Y^{(n)}\to X^{(n)} \leftarrow Z^{(n)})$,
and hence (after unwinding the notation) one sees that the homeomorphism on the right comes from interchanging a colimit and a limit.  To see that the bottom map $\alpha$ is a weak equivalence, note that the maps $Z_n \to X_n$ are Serre fibrations (in each component, this map is just the map from a homotopy orbit space $C_{hU(k)}$ to $BU(k)$), and a colimit of Serre fibrations is a Serre fibration. 
Hence $Z^{(n)}_{\infty} \to X^{(n)}_{\infty}$ (and similarly $Y^{(n)}_{\infty} \to X^{(n)}_{\infty}$) is a Serre fibration.
It is a well-known fact that if $f\co E\to B$ is a Serre fibration, then for each map $g\co A\to B$ there is a weak equivalence
$$\lim \left(A\stackrel{g}{\maps} B \stackrel{f}{\longleftarrow} E \right)
	\stackrel{\heq}{\maps} \holim  \left(A\stackrel{g}{\maps} B \stackrel{f}{\longleftarrow} E \right),
$$
and this precisely tells us that $\alpha$ is a weak equivalence.
Since all of the other maps in Diagram (\ref{tele-holim}) are weak equivalences, so is the top map.

To complete the proof, we must construct filtrations satisfying the four properties listed above. 
Let $C_n(G)\subset \pi_0 \Rep(G)_{hU} = \pi_0 \Rep(G)$ denote the submonoid generated by all representations of dimension at most $n$.  
We define $C_n(G*H)\subset \pi_0 \Rep(G*H)$ to be the pullback of $C_n(G)$ and $C_n(H)$ over
$\bbN = \pi_0 \Rep(\{1\})_{hU}$.  Explicitly, $C_n(G*H)$ is the submonoid consisting of all
components of the form
$$[\oplusdsp_{i=1}^k \rho_i, \oplusdsp_{j=1}^l \psi_j],$$
where all of the representations $\rho_i$ and $\psi_j$ have dimension at most $n$.
We will show that $C_n (G*H)$ is finitely generated.  
By Lemma~\ref{finite-comp}, $\Hom(G, U(n))$ and $\Hom(H, U(n))$ each have finitely many connected components. 
Hence the monoids $C_n(G)$ and $C_n(H)$ have finite generating sets, containing representations of dimension at most $n$.
Choose representatives $\rho_1, \ldots \rho_{r(n)}$ for the components of $\Hom(G, U(m))$, $m = 1, \ldots, n$, and $\psi_1, \ldots \psi_{q(n)}$ for the components of $\Hom(H, U(m))$, $m = 1, \ldots, n$.

We claim that $C_n(G*H)$ is generated by the finite set
$$F_n = \{[\oplusdsp_i \rho_i^{a_i}, \oplusdsp_j \psi_j^{b_j}] \,:\, \mathrm{either\,\,} a_i \leqs n \, \forall i
\mathrm{,\,\, or \,\,} b_j \leqs n \, \forall j\}.$$
(Here we have abbreviated the component associated to the point $[*_{\dim \zeta}, \rho]\in \Rep(\cdot)_{hU}$ by $[\zeta]$.)
Consider a component $C \in C_n (G*H)$.  By definition, 
we may write 
$$C = [\oplusdsp_{l=1}^m \rho_{i_l}, \oplusdsp_{k = 1}^p \psi_{j_k}]$$
where the $\rho_{i_l}$ and the $\psi_{j_k}$ come from our chosen sets of representatives.
We may permute 
the summands of $\oplusdsp_{l=1}^m \rho_{i_l}$ and $\oplusdsp_{k = 1}^p \psi_{j_k}$ independently without changing the component $C$, because $\Hom(G*H, U(n)) \isom \Hom(G, U(n))\cross \Hom(H, U(n))$.  Hence we may write
$$C = [\oplusdsp_{i=1}^{r(n)} \rho_i^{a_i}, \oplusdsp_{j=1}^{q(n)} \psi_j^{b_j}] $$
for some $a_i, b_j\in \bbN$.  If $C\notin F_n$, then there exist $i$ and $j$ such that
$a_i, b_j\geqs n$, and we may assume that $i=j=1$.  Writing
$$C = [\rho_1^{\dim \psi_1}, \psi_1^{\dim \rho_1}] 
	\oplus \left[ \left( \rho_1^{a_1-\dim \psi_1}\right) \oplus \left(\oplusdsp_{i=2}^{r(n)} \rho_i^{a_i}\right),
			\left( \psi_1^{b_1-\dim \rho_1} \right) \oplus \left( \oplusdsp_{j=2}^{q(n)} \psi_j^{b_j} \right)
		\right],
$$
we see that the first factor lies in $F_n$, and by induction on dimension the second factor lies in the submonoid generated by $F_n$.  This completes the proof that $F_n$ generates $C_n (G*H)$.

We now define $\Rep(\cdot)_{hU}^{(n)} = \Rep(\cdot)_{hU}(C_n(\cdot))$.  We must check that these filtrations satisfy properties (1)--(4).  We have explicitly defined our filtrations so that (1) is satisfied.
Property (2) follows from the fact that the monoid $C_n(G*H)$ was defined as a pullback.  
To prove Property (3), note that $\Rep(G*H)_{hU}^{(n)}$ is automatically stably group-like with respect to the sum of (the chosen representatives for) the elements in the generating set $F_n$ (see Example~\ref{fin-gen}).  Letting $[\alpha_n, \beta_n]$ denote this sum, we need to check that $\Rep(G)_{hU}^{(n)}$ is stably group-like with respect to $\alpha_n$ and similarly for $H$.  This follows from the fact the every $\rho_i$ ($i=1,\ldots, r_n$) appears as a summand in $\alpha_n$ (note that some $\psi_j$ must be one-dimensional), and similarly for the $\psi_j$. 
Finally, to check Property (4) we must show that the square
\begin{equation}\label{prop-4}
\xymatrix{
       {\Rep(G * H)_{hU}^{(n)}} \ar[r] \ar[d]
       & {\Rep(G)_{hU}^{(n)}}  \ar[d] \\           
       {\Rep(H)_{hU}^{(n)}} \ar[r]
       & {\Rep(\{1\})_{hU}}  
                    }
\end{equation}
is cartesian.  Since $\Rep(\{1\})_{hU} = \coprod_k BU(k)$, the square (\ref{prop-4}) is the disjoint union of the sub-squares $S_k$ consisting of all points mapping to $BU(k)\subset \Rep(\{1\})_{hU}$.  
Let $C'_n(G)\subset \pi_0 \Rep(G)$ denote the submonoid generated by all representations of dimension at most $n$.  
Let $\Rep(G)^{(n)}_k = \Rep(G)(C'_n(G)) \cap \Hom(G, U(k))$, and define $\Rep(H)^{(n)}_k \subset \Hom(H, U(k))$ similarly.
Let 
\begin{equation*}
\begin{split}
\Rep(G*H)^{(n)}_k = \Rep(G)^{(n)}_k \cross \Rep(H)^{(n)}_k \subset & \,\Rep(G, U(k))\cross \Rep(H, U(k))\\
& = \Rep(G*H, U(k)).
\end{split}
\end{equation*}
Tracing the definitions shows that $S_k$ is in fact the square of homotopy orbit spaces associated
to the cartesian square
\begin{equation*}
\xymatrix{
       {\Rep(G * H)^{(n)}_k} \ar[r] \ar[d]
       & {\Rep(G)^{(n)}_k}  \ar[d] \\           
       {\Rep(H)^{(n)}_k} \ar[r]
       & {*,}  
                    }
\end{equation*}
so Lemma~\ref{cartesian-orbit} implies that $S_k$ is cartesian.  Hence the square (\ref{prop-4}) is a disjoint union of cartesian squares.  This completes the proof of Theorem~\ref{excision-free-prod}.
\end{proof}

As an application of Theorem~\ref{excision-free-prod}, we now compute $\K^*(PSL_2(\bbZ))$ in the unitary case.  It is well known that $PSL_2(\bbZ) \isom \bbZ/2 * \bbZ/3$ (see Alperin~\cite{Alperin} for a short proof).  We work with $PSL_2 (\bbZ)$ for concreteness, although the same argument gives a computation of $\K^*(G*H)$ for all finite groups $G$ and $H$.

T. Lawson~\cite[Chapter 6.3]{Lawson-thesis} has shown, using basic representation theory, that for every finite group $G$, $\K(G)\heq \bigvee_k \mathbf{ku}$, where $k$ is the number of irreducible representations of $G$.  Hence in particular 
$$\K^*(\bbZ/m) =  \left\{ \begin{array}{ll}
					 0, \,\,\, * \mathrm{\,\, odd,}\\
					 \bbZ^m, \,\,\, * \mathrm{\,\, even.}
				      \end{array}
					\right.
$$
By Remark~\ref{cart-htpy-gps}, we now have (for each $i\geqs 0$) $\K^{2i+1} (\bbZ/2 * \bbZ/3) \isom 0\oplus 0 = 0$ and an exact sequence
$$0\maps \K^{2i} (\bbZ/2* \bbZ/3) \maps \bbZ^2 \oplus \bbZ^3 \maps \bbZ \maps 0,$$
from which it follows that $\K^{2i}(\bbZ/2* \bbZ/3) \isom \bbZ^4$.  Thus we have:

\begin{proposition}$\label{PSL}$
$$\K^*(PSL_2(\bbZ)) =  \left\{ \begin{array}{ll}
					 0, \,\,\, *  \mathrm{\,\, odd,}\\
					 \bbZ^4, \,\,\, * \mathrm{\,\, even.}
				      \end{array}
					\right.
$$
\end{proposition}

We briefly indicate Lawson's computation of $\K(G)$ for $G$ finite.  Every representation $\rho$ of $G$ breaks up canonically into isotypical components, and together with Schur's Lemma this gives a permutative functor $\mcR(G)\to \Vect^k$, which records the dimensions of the isotypical components.  Here $k$ is the number of irreducible representations of $G$ and $\Vect$ is the category with $\N$ as objects and $\coprod_n U(n)$ as morphisms (we will work in the unitary case, but the general linear case is identical).  This functor is continuous, since any two representations connected by a path are isomorphic (since $G$ is finite, the trace gives a continuous, complete invariant of the isomorphism type, and it can take on only countably many values).  One now checks that this functor induces a weak equivalence on classifying spaces, and hence on $K$--theory spectra.  This is rather like the proof of Proposition~\ref{classifying-space}: one sees that $B\mcR(G)$ is a model for $\coprod_{\rho_i} B \left( \Stab (\rho_i) \right)$, where the $\rho_i$ are representatives for the isomorphism types.  Now Schur's Lemma implies that $\Stab (\rho_i) \isom \prod U(n_j)$, where the $n_j$ are the dimensions of the isotypical components of $\rho_i$.  The comparison with
$B(\Vect^k) \homeo B(\Vect)^k \homeo \left(\coprod_n BU(n) \right)^k$ is now straightforward.

%
%
%
\bibliographystyle{gtart}

\bibliography{references}

\end{document}